# Probability of Failure of Safety-Critical Systems Subject to Partial Tests

Florent Brissaud, Institut National de l'Environnement Industriel et des Risques
Anne Barros, Université de Technologie de Troyes
Christophe Bérenguer, Université de Technologie de Troyes

*SUMMARY & CONCLUSIONS*

A set of general formulas is proposed for the probability of failure on demand (PFD) assessment of MooN architecture (i.e. k-out-of-n) systems subject to proof tests. The proof tests can be partial or full. The partial tests (e.g. visual inspections, partial stroke testing) are able to detect only some system failures and leave the others latent, whereas the full tests refer to overhauls which restore the system to an as good as new condition. Partial tests may occur at different time instants (periodic or not), up to the full test. The system performances which are investigated are the system availability according to time, the PFD average in each partial test time interval, and the total PFD average calculated on the full test time interval.

Following the given expressions, parameter estimations are proposed to assess the system failure rates and the partial test effectiveness according to feedback data from previous test policies. Subsequently, an optimization of the partial test strategy is presented. In the 2oo6 system given as example, an improvement of about 10% of the total PFD average has been obtained, just by a better (non-periodic) distribution of the same number of partial tests, in the full test time interval.

## 1 INTRODUCTION

The safety instrumented systems (SIS) play a major part in industrial risk management as safety barriers. A SIS aims at performing at least one safety function to achieve or maintain a safe state of equipment under control (EUC), in respect of a specific hazardous event. Due to the critical role of SIS for health, environment, and goods, functional safety standards have been developed as, for example, the IEC 61508 [1] and the ANSI/ISA S84.00.01-2004 [2]. A framework is then introduced to consider the overall system and software safety lifecycle and some requirements aim at avoiding the systematic failures. In the realization phase, the safety integrity level (SIL) is determined with respect to the random hardware failures. For a SIS operating in a low demand mode (i.e. the safety function is only performed on demand, and not more than once per year), the used criterion is the average probability of system failure to perform its safety function on demand. It corresponds to the mean system unavailability, calculated on the full test time interval, and denoted PFDavg.

According to IEC 61508 [1] and to the assumptions given afterwards, different characteristics have to be taken into account for PFDavg assessment: system architecture, failure rates, proof test intervals and effectiveness. In the present paper, a proof test is described as "partial" if it is imperfect, that is, it is able to detect only some system failures and leave the others latent. Visual inspections and partial stroke testing are examples of such partial tests. When a proof test is perfect (i.e. its effectiveness is equal to one), it is described as "full" and refers to overhauls which restore the system to an as good as new condition. Even if partial tests are less effective, they can be preferred to the full tests for several reasons:

- Full tests are generally physical (e.g. stimulation of sensing elements), costly and time consuming. They may sometimes be substituted by electronic tests (e.g. electronic simulation) but do not cover all the failures.
- Full tests often imply stopping production (e.g. power supply cut, flow stop by safety valves), sometimes unacceptable for industrialists. Partial tests (e.g. quarter-turn valve closure) are therefore preferred.
- Some safety devices cannot be fully tested without degradations or destructions (e.g. one shot devices).
- Only real conditions testing can pretend to be full and, in many cases, may provoke more hazards than prevention (e.g. fire, toxic gas, or overpressure detection).

Methods for PFDavg assessment are mentioned in [1] (e.g. fault trees, reliability block diagrams, Markov models) but none of them is prescribed. These techniques are compared and discussed in [3]-[5]. For example, a Markov model is used in [6] and reliability block diagrams in [7], but both of them ignore the partial tests. The use of such tests can be regarded as a kind of imperfect preventive maintenance. Various methods and optimal policies for imperfect maintenance have been discussed in [8]. For example, some approaches set a constant probability that a preventive maintenance is perfect or imperfect [9]-[10]. More specifically, a probability that a failure remained undiscovered after testing is assumed in [11], but only for a single unit. References [3]-[4] claim that enhanced Markov analyses cover most aspects relevant for safety quantification. The use of Markov models is also defended in [5]. Because partial and full tests usually occur at deterministic time instants (e.g. periodically), basic Markov models are not appropriate [12] and extended Markov models have been developed [12]-[14]. Similar approaches have been used for cost optimization [15].

In the present paper, a set of general formulas, given in a neat form and easy to compute, has been obtained by analytical approaches. These expressions are proposed to assess the availability of MooN architecture (i.e. k-out-of-n) systems made up by homogeneous components (i.e. components with identical failure rates), and subject to partial and full tests. The partial tests may occur at different time instants (periodic or not), up to the full test. Section 2 presents the assumptions, notations, and formulas. An application is then given in Section 3 for parameter estimations, PFDavg assessment, and optimization of the partial test distribution.

## 2 PFD ASSESSMENT

### 2.1 General Assumptions

- All failures taken into account are dangerous and only detected by partial or full tests.
- The system is made up by $N$ components which are independent and have identical and constant failure rates.
- The system has a MooN architecture i.e. it is made up by $N$ components and it is able to perform its safety function if any $M$ or more components among $N$ are in an operating state (i.e. $M$-out-of-$N$ system).
- The $N$ components of the system are in an operating state at time $t_0$.
- The partial tests are able to detect only some specific failures of each component of the system.
- The full tests are able to detect all the failures of each component of the system.
- All components are tested together during any test.
- When a failure is detected by a partial or a full test, it is repaired immediately. During test and maintenance actions, measures are performed to maintain the EUC in a safe state in such a way that any test or maintenance duration is not included in the analyses.
- After each full test, the system is restored to an as good as new condition. The PFDavg can therefore be assessed according to the full test time interval.

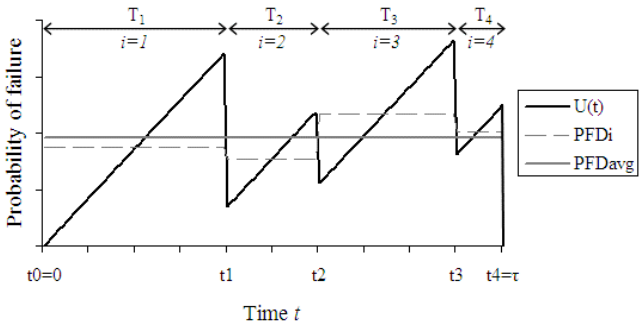

Figure 1 – Notation, example with n = 4

### 2.2 Notations

| | |
|---|---|
| MooN | system architecture, with $M \leq N$ |
| $\lambda$ | failure rate of any of the $N$ components which compose the system |
| $A_e(t)$ | availability function of any of the $N$ components which compose the system i.e. probability that the component is in an operating state at time $t$ |
| $A(t)$ | availability function of the system i.e. probability that the system is able to perform its safety function at time $t$ |
| $U(t)$ | unavailability function of the system i.e. $U(t)=1-A(t)$ |
| $t_i$ | time instant of the $i^{th}$ test (which can be partial or full), with the initial condition $t_0=0$ |
| $T_i$ | time interval between the $(i-1)^{th}$ and the $i^{th}$ test i.e. $T_i=t_i-t_{i-1}$ |
| $E$ | efficiency of partial tests i.e. a proportion equal to $E$ of each component failure rate corresponds to failures which are detected by any partial test |
| $n$ | total number of tests in the full test time interval i.e. $(n-1)$ partial tests plus the $n^{th}$ test which is full |
| $\tau$ | full test time interval i.e. $\tau=t_n$ |
| $PFD_i$ | average probability of the system failure to perform its safety function on demand (i.e. mean unavailability) in the time interval between the $(i-1)^{th}$ and the $i^{th}$ test (i.e. $[t_{i-1}, t_i]$) |
| $PFD_{avg}$ | average probability of the system failure to perform its safety function on demand (i.e. mean unavailability) in the full test time interval (i.e. $[0, \tau]$) |
| $Obs_i$ | probability, for each component, to observe a failure during the $i^{th}$ test |
| $k_i$ | number of components observed in a failed state during the $i^{th}$ test |
| $K$ | equivalent total number of components observed during each test |

Notations are reported in Figure 1. A system is defined by the set $\{M, N, \lambda\}$, and a test policy can be either defined by the set $\{E, t_1, t_2, ..., t_n\}$ or by the set $\{E, T_1, T_2, ..., T_n\}$.

### 2.3 General Formulas

For each component, a part with a failure rate equal to $E \cdot \lambda$ is testable by any partial or full test, and another part with failure rate equal to $(1-E) \cdot \lambda$ is only testable by the full tests. The corresponding reliability block diagram (RBD) is then given in Figure 2, and (see proof in Appendix):

$$A_e(t) = e^{E \cdot \lambda \cdot t_{i-1}} \cdot e^{-\lambda \cdot t} \quad \text{for } t \in [t_{i-1}, t_i[ \quad (1)$$

The availability function of the system is therefore, (see proof in Appendix):

$$A(t) = \sum_{x=M}^{N} \left[ S(M,N,x) \cdot e^{x \cdot E \cdot \lambda \cdot t_{i-1}} \cdot e^{-x \cdot \lambda \cdot t} \right] \quad \text{for } t \in [t_{i-1}, t_i[ \quad (2)$$

With the following time-independent sum (see Table I):

$$S(M,N,x) = \sum_{k=M}^{x} \left[ \binom{N}{x} \cdot \binom{x}{k} \cdot (-1)^{x-k} \right]$$

$$\text{for } x = M,...,N \quad (3)$$

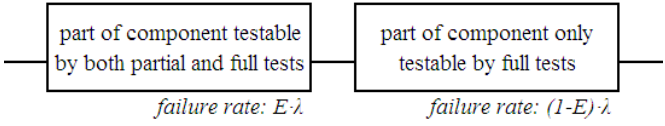

*Figure 2 – Reliability Block Diagram (RBD) for any of the N components which compose the system*

| | $N=1$ | $N=2$ | $N=3$ | $N=4$ |
|---|---|---|---|---|
| $M=1$ | 1 | 2 / -1 | 3 / -3 / 1 | 4 / -6 / 4 / -1 |
| $M=2$ | | 1 | 3 / -2 | 6 / -8 / 3 |
| $M=3$ | | | 1 | 4 / -3 |
| $M=4$ | | | | 1 |

*Table 1 – S(M,N,M) / S(M,N,M+1) / ... / S(M,N,N) values for some MooN architectures*

The average probabilities of the system failure to perform its safety function on demand, in the time interval between the $(i-1)^{th}$ and the $i^{th}$ test, and in the full test time interval, are therefore (see proof in Appendix):

$$PFD_i = 1 - \sum_{x=M}^{N} \left[ S(M,N,x) \cdot e^{-x \cdot (1-E) \lambda \cdot t_{i-1}} \cdot \frac{1 - e^{-x \cdot \lambda \cdot T_i}}{x \cdot \lambda \cdot T_i} \right] \quad (4)$$

$$PFD_{avg} = 1 - \sum_{x=M}^{N} \left[ S(M,N,x) \cdot \sum_{i=1}^{n} \left[ e^{-x \cdot (1-E) \lambda \cdot t_{i-1}} \cdot \frac{1 - e^{-x \cdot \lambda \cdot T_i}}{x \cdot \lambda \cdot \tau} \right] \right] \quad (5)$$

When $\lambda \cdot \tau$ is small (i.e. $\lambda \cdot \tau << 10^{-2}$), the following approximations can be done, using the Taylor's theorem:

$$A_e(t) \approx 1 + E \cdot \lambda \cdot t_{i-1} - \lambda \cdot t \quad \text{for } t \in [t_{i-1}, t_i[ \quad (6)$$

$$A(t) \approx 1 - \binom{N}{M-1} \cdot \lambda^{N-M+1} \cdot (t - E \cdot t_{i-1})^{N-M+1}$$
$$\text{for } t \in [t_{i-1}, t_i[ \quad (7)$$

$$PFD_i \approx \binom{N}{M-1} \cdot \frac{\lambda^{N-M+1}}{N-M+2} \cdot \frac{1}{T_i} \quad (8)$$
$$\cdot \left( (t_i - E \cdot t_{i-1})^{N-M+2} - (t_{i-1} \cdot (1-E))^{N-M+2} \right)$$

$$PFD_{avg} \approx \binom{N}{M-1} \cdot \frac{\lambda^{N-M+1}}{N-M+2} \cdot \frac{1}{\tau}$$
$$\cdot \sum_{i=1}^{n} \left[ (t_i - E \cdot t_{i-1})^{N-M+2} - (t_{i-1} \cdot (1-E))^{N-M+2} \right] \quad (9)$$

### 2.4 Special Case: Without Partial Test

Without partial test, (2), (7), (5), and (9) become:

$$A(t)^{(w)} = \sum_{x=M}^{N} S(M,N,x) \cdot e^{-x \cdot \lambda \cdot t} \quad \text{for } t \in [0, \tau[ \quad (10)$$

$$A(t)^{(w)} \approx 1 - \binom{N}{M-1} \cdot (\lambda \cdot t)^{N-M+1} \quad \text{for } t \in [0, \tau[ \quad (11)$$

$$PFD_{avg}^{(w)} = 1 - \sum_{x=M}^{N} \left[ S(M,N,x) \cdot \frac{1 - e^{-x \cdot \lambda \cdot \tau}}{x \cdot \lambda \cdot \tau} \right] \quad (12)$$

$$PFD_{avg}^{(w)} \approx \binom{N}{M-1} \cdot \frac{(\lambda \cdot \tau)^{N-M+1}}{N-M+2} \quad (13)$$

### 2.5 Special Case: With Periodic Partial Tests

Partial tests are periodic if $T_i = T_0$ for $i=1,...,n$, then $t_i = i \cdot T_0$ for $i=1,...,n$. Formulas (2), (7), (5), and (9) therefore become:

$$A(t)^{(p)} = \sum_{x=M}^{N} S(M,N,x) \cdot e^{x \cdot E \cdot \lambda \cdot (i-1) T_0} \cdot e^{-x \cdot \lambda \cdot t}$$
$$\text{for } t \in [(i-1) \cdot T_0, i \cdot T_0[ \quad (14)$$

$$A(t)^{(p)} \approx 1 - \binom{N}{M-1} \cdot \lambda^{N-M+1} \cdot (t - E \cdot (i-1) \cdot T_0)^{N-M+1}$$
$$\text{for } t \in [(i-1) \cdot T_0, i \cdot T_0[ \quad (15)$$

$$PFD_{avg}^{(p)} = 1 - \sum_{x=M}^{N} \left[ S(M,N,x) \cdot \frac{1 - e^{-x \cdot \lambda \cdot T_0}}{x \cdot \lambda \cdot T_0} \cdot \frac{1}{n} \cdot \sum_{i=1}^{n} \left[ e^{-x \cdot (1-E) \lambda \cdot (i-1) T_0} \right] \right] \quad (16)$$

$$PFD_{avg}^{(p)} \approx \binom{N}{M-1} \cdot \frac{(\lambda \cdot T_0)^{N-M+1}}{N-M+2} \cdot \frac{1}{n}$$
$$\cdot \sum_{j=0}^{n-1} \left[ (1 + j \cdot (1-E))^{N-M+2} - (j \cdot (1-E))^{N-M+2} \right] \quad (17)$$

## 3 APPLICATION

### 3.1 Case Study

A system for oxygen concentration measurement is used as case study. It takes part of an inerting system which aims to reduce the oxygen level in atmospheric air by introducing nitrogen in controlled amounts. The oxygen concentration is maintained below a high level which does not allow the fire breaking out, and above a low level which remains the place accessible to people. To control the oxygen concentration, six oxygen transmitters are used. Because the nitrogen is quickly and heterogeneously distributed into the air space, these six transmitters are assumed redundant. The safety function then consists in detecting a low or a high oxygen level, according to a 2oo6 architecture. Two transmitter test procedures are recommended by the manufacturer:

- Every year: control of the measurements and, if required, followed by adjustment i.e. full test.
- Occasionally: visual inspections and some electronic checks i.e. partial test.

For cost reasons, all the six transmitters are tested together during any test. The basic test policy consists of a full test every year and a partial test every three months.

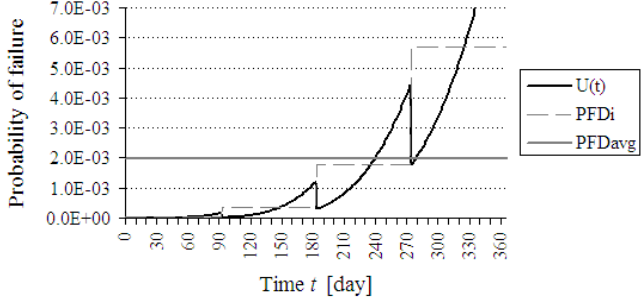

*Figure 3 – PFD for the case study, according to the basic test policy*

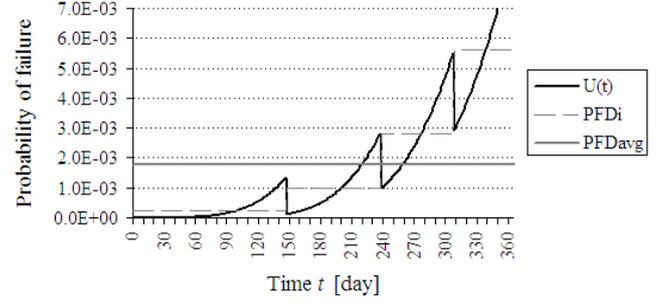

*Figure 4 – PFD for the case study, according to the optimized test policy*

The first step of the application uses feedback data from the partial and full tests, in order to estimate the transmitter failure rates, and the partial test effectiveness. Then, the probabilities of system failure are assessed according to this basic test policy. Subsequently, the partial test distribution is optimized in order to reduce the PFDavg.

### 3.2 Parameters Estimation and PFD Assessment

During the $i^{th}$ test of each transmitter, the probability to observe a failure is, according to the RBD given in Figure 2:

$$Obs_i \approx E \cdot \lambda \cdot T_i \quad \text{for } i = 1,\ldots,(n-1) \quad (18)$$

$$Obs_i \approx E \cdot \lambda \cdot T_i + (1-E) \cdot \lambda \cdot \tau \quad \text{for } i = n \quad (19)$$

Moreover, because a number of components $k_i$ among $K$ have been observed in a failed state during the $i^{th}$ test, an intuitive and empirical estimator of $Obs_i$ is:

$$\hat{Obs}_i = k_i / K \quad \text{for } i = 1,\ldots,n \quad (20)$$

The following estimators for the failure rate $\lambda$ and the partial test effectiveness $E$ can be deduced from (18)-(20):

$$\hat{\lambda} = \frac{1}{K \cdot \tau} \cdot \left( \sum_{i=1}^{n} k_i \right) \quad (21)$$

$$\hat{E} = \frac{\tau}{t_{n-1}} \cdot \left( \sum_{i=1}^{n-1} k_i \Big/ \sum_{i=1}^{n} k_i \right) \quad (22)$$

The observations $k_i$ follow a binomial distribution. Confidence intervals for the estimations given by (21) and (22) can therefore be obtained using Fisher distributions [16].

Assuming four installations, each of them using six oxygen transmitters, for a time period equal to four years, the equivalent total number of transmitters observed during each test is thus $K=4 \cdot 6 \cdot 4=96$. A total number of $k_1+k_2+k_3=16$ transmitter failures have been observed during the partial tests, and $k_4=35$ transmitter failures during the full tests. According to (21) and (22), the transmitter failure rate $\lambda$ is estimated at $6.1 \cdot 10^{-5}$ hour$^{-1}$, and the partial test effectiveness $E$ at $0.42$.

According to (16), equivalent to (5) for periodic partial tests, the PFDavg is then equal to $2.06 \cdot 10^{-3}$. For the present case study, the partial tests have reduced the PFDavg by a factor of five, compared to the result without partial test given by (12). The probabilities of system failure, according to this basic test policy, are depicted in Figure 3.

### 3.3 Optimization of the Partial Tests Distribution

An optimization of the test policy consists in distributing the same number of partial tests, inside the full test time interval, in order to minimize the PFDavg. When the costs related to the partial tests are independent of the time instants, this approach may improve the system safety without additional cost, provided that the test times are optimally chosen.

The optimal partial test time instants are denoted $t_i^*$ with $i=1,\ldots,(n-1)$, and the optimal partial test time intervals are denoted $T_i^*$ with $T_i^*=t_i^*-t_{i-1}^*$. The optimized test policy is then obtained by solving the following equation:

$$\{t_1^*, t_2^*, \ldots, t_{n-1}^*\} = \arg \min_{t_1 \leq t_2 \leq \ldots \leq t_{n-1}} \left[ PFD_{avg} \right] \quad (23)$$

With the PFDavg as defined by (5).

For the present case study, the optimal partial test time instants and the corresponding partial test time intervals, obtained by solving (23), are reported in Table 2 (last row).

Using this optimized test policy, and according to (5), the PFDavg is now equal to $1.87 \cdot 10^{-3}$, that is a reduction of about 10% compared to the result with the basic (periodic) test policy. Moreover, further analyses show that the maximum system unavailability in the full test time interval has been reduced by more than 25%. The probabilities of system failure, according to this optimized test policy, are depicted in Figure 4.

| Test Policy | Partial Test Distribution | | PFDavg |
|---|---|---|---|
| Without partial test | $\tau=12.0$ | | $1.03 \cdot 10^{-2}$ |
| Basic test policy (periodic) | $t_1=3.0$<br>$t_2=6.0$<br>$t_3=9.0$<br>$t_4=\tau=12.0$ | $T_1=3.0$<br>$T_2=3.0$<br>$T_3=3.0$<br>$T_4=3.0$ | $2.06 \cdot 10^{-3}$ |
| Optimized test policy | $t_1^*=4.8$<br>$t_2^*=7.8$<br>$t_3^*=10.1$<br>$t_4^*=\tau=12.0$ | $T_1^*=4.8$<br>$T_2^*=3.0$<br>$T_3^*=2.3$<br>$T_4^*=1.9$ | $1.87 \cdot 10^{-3}$ |

*Table 2 – Results summary*

## 4 CONCLUSION

By introducing general formulas for the probability of failure on demand (PFD) assessments of MooN architecture systems subject to partial and full tests, the proposed work provides practical tools for risk management. The neat form of the proposed exact and approximate expressions allows the system and test policy performances to be estimated and optimized quite simply and directly. It has been particularly shown that the average probability of system failure to perform its safety function on demand (PFDavg) can be reduced, just by a better (non-periodic) distribution of an appointed number of partial tests. Thus offering good prospects are for safety improvement without additional cost.

Most of the general assumptions given in Section 2.1 are easily fulfilled when measures are performed to maintain the EUC in a safe state during test and maintenance actions. However, some other aspects should be generalized in order to be applicable to a wider field of applications, as for example:
- Systems with heterogeneous components (i.e. components with different failure rates).
- Systems subject to common cause failures.
- Systems with extended MooN architectures (e.g. MooN architecture parts into each other, capacitated systems).
- Aging systems.
- Stochastic test policies, etc.

Most of these previous cases can be solved using a similar reasoning as proposed in this paper (see Appendix), but make the general formulas sometimes much more difficult to grasp. Another prospect for development concerns the staggered tests which allow for reducing probabilities of system failure by testing the system components at different time instants [17].

## APPENDIX

The following results are for $t_{i-1} \leq t < t_i$ with $i=1, ..., n$.

- *Proof of (1)*: according to the RBD given in Figure 2:

$$A_e(t) = e^{-E \cdot \lambda \cdot (t-t_{i-1})} \cdot e^{-(1-E)\lambda \cdot t} = e^{E \cdot \lambda \cdot t_{i-1}} \cdot e^{-\lambda \cdot t} \quad (24)$$

- *Proof of (2)*: according to [16]:

$$A(t) = \sum_{k=M}^{N} \left[ \binom{N}{k} \cdot A_e(t)^k \cdot \left(1 - A_e(t)\right)^{N-k} \right] \quad (25)$$

$$A(t) = \sum_{k=M}^{N} \left[ \binom{N}{k} \cdot e^{k \cdot E \cdot \lambda \cdot t_{i-1}} \cdot e^{-k \cdot \lambda \cdot t} \cdot \left(1 - e^{E \cdot \lambda \cdot t_{i-1}} \cdot e^{-\lambda \cdot t}\right)^{N-k} \right] \quad (26)$$

Using the Newton's binomial theorem:

$$A(t) = \sum_{k=M}^{N} \left[ \binom{N}{k} \cdot e^{k \cdot E \cdot \lambda \cdot t_{i-1}} \cdot e^{-k \cdot \lambda \cdot t} \cdot \sum_{l=0}^{N-k} \left[ \binom{N-k}{l} \cdot (-1)^{N-k-l} \cdot \left(e^{(N-k-l)E \cdot \lambda \cdot t_{i-1}} \cdot e^{-(N-k-l)\lambda \cdot t}\right) \right] \right] \quad (27)$$

$$A(t) = \sum_{k=M}^{N} \sum_{l=0}^{N-k} \left[ \binom{N}{k} \cdot \binom{N-k}{l} \cdot (-1)^{N-k-l} \cdot e^{(N-l)E \cdot \lambda \cdot t_{i-1}} \cdot e^{-(N-l)\lambda \cdot t} \right] \quad (28)$$

Using the Fubini's theorem:

$$A(t) = \sum_{l=0}^{N-M} \sum_{k=M}^{N-l} \left[ \binom{N}{k} \cdot \binom{N-k}{l} \cdot (-1)^{N-k-l} \cdot e^{(N-l)E \cdot \lambda \cdot t_{i-1}} \cdot e^{-(N-l)\lambda \cdot t} \right] \quad (29)$$

By substituting $x=N-l$:

$$A(t) = \sum_{x=M}^{N} \sum_{k=M}^{x} \left[ \binom{N}{x} \cdot \binom{x}{k} \cdot (-1)^{x-k} \cdot e^{x \cdot E \cdot \lambda \cdot t_{i-1}} \cdot e^{-x \cdot \lambda \cdot t} \right] \quad (30)$$

Finally:

$$A(t) = \sum_{x=M}^{N} \left[ S(M,N,x) \cdot e^{x \cdot E \cdot \lambda \cdot t_{i-1}} \cdot e^{-x \cdot \lambda \cdot t} \right] \quad (31)$$

With:

$$S(M,N,x) = \sum_{k=M}^{x} \left[ \binom{N}{x} \cdot \binom{x}{k} \cdot (-1)^{x-k} \right] \quad (32)$$

- *Proofs of (4)*:

$$PFD_i = \frac{1}{T_i} \cdot \int_{t_{i-1}}^{t_i} U(t) \cdot dt = 1 - \frac{1}{T_i} \cdot \int_{t_{i-1}}^{t_i} A(t) \cdot dt \quad (33)$$

$$PFD_i = 1 - \sum_{x=M}^{N} \left[ S(M,N,x) \cdot e^{-x \cdot (1-E)\lambda \cdot t_{i-1}} \cdot \frac{1 - e^{-x \cdot \lambda \cdot T_i}}{x \cdot \lambda \cdot T_i} \right] \quad (34)$$

- *Proof of (5)*:

$$PFD_{avg} = \frac{1}{\tau} \cdot \sum_{i=1}^{n} \left[ T_i \cdot PFD_i \right] \quad (35)$$

$$PFD_{avg} = 1 - \sum_{x=M}^{N} \left[ S(M,N,x) \cdot \sum_{i=1}^{n} \left[ e^{-x \cdot (1-E)\lambda \cdot t_{i-1}} \cdot \frac{1 - e^{-x \cdot \lambda \cdot T_i}}{x \cdot \lambda \cdot \tau} \right] \right] \quad (36)$$

## REFERENCES

1. International Electrotechnical Commission, *IEC 61508 Functional safety of electrical / electronic / programmable electronic safety-related systems*, Geneva, IEC Standard, 2002.
2. International Society of Automation, *ANSI/ISA-84.00.01-2004 Functional safety: Safety Instrumented Systems for the Process Industry Sector*, Research Triangle Park, ISA Standard, 2004.
3. J. L. Rouvroye, A. C. Brombacher, "New quantitative safety standards: different techniques, different results?" *Reliability Engineering & System Safety*, vol. 66, 1999, pp. 121-125.

4. J. L. Rouvroye, E. G. van den Bliek, "Comparing safety analysis techniques," *Reliability Engineering & System Safety*, vol. 75, 2002, pp. 289-294.
5. J. V. Bukowski, "A comparison of techniques for computing PFD average," *Proc. Ann. Reliability & Maintenability Symp.*, 2005, pp. 590-595.
6. T. Zhang, W. Long, Y. Sato, "Availability of systems with self-diagnosis components–applying Markov model to IEC 61508-6," *Reliability Engineering & System Safety*, vol. 80, 2003, pp. 133-141.
7. H. Guo, X. Yang, "A simple reliability block diagram method for safety integrity verification," *Reliability Engineering & System Safety*, vol. 92, 2007, pp. 1267-1273.
8. H. Pham, H. Wang, "Imperfect maintenance," *European Journal of Operational Research*, vol. 94, 1996, pp. 425-438.
9. M. Brown, F. Porschan, "Imperfect repair," *Journal of Applied Probability*, vol. 20, 1983, pp. 851-859.
10. T. Nakagawa, K. Yasui, "Optimal policies for a system with imperfect maintenance," *IEEE Trans. Reliability*, vol. 36, 1987, pp. 631-633.
11. F. G. Badía, M. D. Berrade, C. A. Campos, "Optimal inspection and preventive maintenance of units with revealed and unrevealed failures," *Reliability Engineering & System Safety*, vol. 78, 2002, pp. 157-163.
12. J. V. Bukowski, "Modeling and analyzing the effects of periodic inspection on the performance of safety-critical systems," *IEEE Trans. Reliability*, vol. 50, 2001, pp. 321-329.
13. G. Levitin, T. Zhang, M. Xie "State probability of a series-parallel repairable system with two-types of failure states," *International Journal of Systems Sciences*, vol. 37, 2006, pp. 1011-1020.
14. M. Kumar, A. Verma, A. Srividya, "Modeling demand rate and imperfect proof-test and analysis of their effect on system safety," *Reliability Engineering & System Safety*, vol. 93, 2008, pp. 1720-1729.
15. C. Wang, S. H. Sheu, "Determining the optimal production maintenance policy with inspection errors: using a Markov chain," *Computers & Operations Research*, vol. 30, 2003, pp. 1-17.
16. M. Rausand, A. Høyland, *System Reliability Theory: Models, Statistical Methods, and applications – Second Edition*, New Jersey, Wiley & Sons, 2004.
17. J. L. Rouvroye, J. A. Wiegerinck, "Minimizing costs while meeting safety requirements: Modeling deterministic (imperfect) staggered tests using standard Markov models for SIL calculations," *ISA Trans.*, vol. 45, 2006, pp. 611-621.

*BIOGRAPHIES*

Florent Brissaud
Institut National de l'Environnement Industriel et des Risques
Parc Technologique ALATA – BP 2
60550 Verneuil-en-Halatte, France

e-mail: florent.brissaud@ineris.fr

Florent Brissaud received his MS in Dependability, Risk Analysis, and Environment at the Troyes University of Technology (UTT), in France. Florent is currently a Safety & Reliability Engineer at the French National Institute for Industrial Environment and Risk (INERIS), in the Safety Barrier Evaluation Unit (BT2S). He started his PhD thesis in 2007 under the scientific supervision of the UTT. His research interests mainly focus on reliability evaluation of safety instrumented systems based on new technologies.

Anne Barros, PhD
ICD, FRE CNRS 2848
Université de Technologie de Troyes
12 rue Marie Curie – BP 2060
F-10010 Troyes cedex, France

e-mail: anne.barros@utt.fr

Anne Barros is Associate Professor at the Troyes University of Technology (UTT) since 2004 and is member of the Systems Modeling and Dependability Laboratory. Her research interests focus on the construction of probabilistic decision indicators used for safety assessment and maintenance optimization. These indicators are based on stochastic modeling of the system state, its environment, and generally speaking of the whole information given by the monitoring device. She has authored papers on stochastic models in reliability or maintenance optimization in IEEE transactions on reliability, Reliability Engineering and System Safety, Journal of Risk and Reliability, Journal of Loss Prevention in the Process Industries, International Journal of reliability, Quality and Safety Engineering.

Christophe Bérenguer, PhD
ICD, FRE CNRS 2848
Université de Technologie de Troyes
12 rue Marie Curie – BP 2060
F-10010 Troyes cedex, France

e-mail: christophe.berenguer@utt.fr

Christophe Bérenguer is Professor and Head of the Industrial Engineering Program and of the System Safety and Optimization PhD program at the Université de Technologie de Troyes (France). He has served as an officer (treasurer) for the European Safety and Reliability Association since 2005. He is member of the Editorial Boards of Reliability Engineering and System Safety and Journal of Risk and Reliability. His research interests include stochastic modeling of systems deterioration, performance evaluation and optimization of maintenance policies, health monitoring and probabilistic safety assessment.